\input amstex
\documentstyle{amsppt}
\nopagenumbers
\accentedsymbol\tx{\tilde x}
\accentedsymbol\ty{\tilde y}
\def\msum{\operatornamewithlimits{\sum\!{\ssize\ldots}\!\sum}}
\def\tr{\operatorname{tr}}
\def\Res{\operatornamewithlimits{Res}}
\pagewidth{13cm}
\pageheight{21cm}
\rightheadtext{On some equations that can be brought \dots}
\topmatter
\title On some equations that can be brought to
the equations of diffusion type.
\endtitle
\author
Dmitrieva~V.~V., Gladkov~A.~V., Sharipov~R.~A.
\endauthor
\abstract For the system of second order quasilinear parabolic 
equations the problem of reducing them to the equations 
of diffusion type is considered. In non-degenerate case an 
effective algorithm for solving this problem is suggested.
\endabstract
\thanks Authors are grateful to International Soros Foundation 
(``Open Society'' institute), to Russian Fund for Fundamental 
Research, and the Academy of Sciences of the Republic Bashkortostan
for financial support in 1998.
\endthanks
\address
Department of Mathematics, Ufa State Aviation Technical University,
Karl Marks str\. 12, 450000 Ufa, Russia.
\endaddress
\email
gladkov\@math.ugatu.ac.ru
\endemail
\address
Department of Mathematics, Bashkir State University,
Frunze str\. 32, 450074 Ufa, Russia.
\endaddress
\email \vtop{\hsize 6cm \noindent DmitrievaVV\@ic.bashedu.ru\newline
R\_\hskip 0.7pt Sharipov\@ic.bashedu.ru}
\endemail
\endtopmatter
\loadbold
\document
\head
1. Introduction.
\endhead
     Let's consider the system of differential equations defined
by some matrix $A=A(y^1,\ldots,y^n)$ with non-zero determinant:
$$
\frac{\partial y^i}{\partial\,t}=\sum^n_{j=1}\frac{\partial}{\partial x}
\left(A^i_j\,\frac{\partial y^j}{\partial x}
\right),\text{\ \ where \ }i=1,\ldots,n.\hskip -2em
\tag1.1
$$
The equations \thetag{1.1} do not form the set invariant under the 
point transformations given by the following changes of variables:
$$
\align
\ty^1=\ty^1(y^1,\ldots,y^n),\hskip -2em\\
.\ .\ .\ .\ .\ .\ .\ .\ .\ .\ .\ .\ .\ .\ \hskip -2em
\tag1.2\\
\ty^n=\ty^1(y^1,\ldots,y^n).\hskip -2em
\endalign
$$
However, one can expand this set up to the point-invariant class
of equations: 
$$
\frac{\partial y^i}{\partial\,t}=\sum^n_{j=1}A^i_j\left(
\frac{\partial^2 y^j}{\partial x^2}+\shave{\sum^n_{r=1}\sum^n_{s=1}}
\Gamma^j_{rs}\,\frac{\partial y^r}{\partial x}\,\frac{\partial y^s}
{\partial x}\right).\hskip -2em
\tag1.3
$$
Here coefficients $A^i_j$ and $\Gamma^j_{rs}$ depend on $y^1,\,\ldots,
\,y^n$ so that $\det A\neq 0$. Denote by $T$ the Jacoby matrix for
the change of variables \thetag{1.2} and denote by $S$ the Jacoby matrix
for inverse change of variables:
$$
\xalignat 2
&S^i_j=\frac{\partial y^i}{\partial\ty^j},
&&T^i_j=\frac{\partial\ty^i}{\partial y^j}.
\endxalignat
$$
In geometry matrices $S$ and $T$ are called direct and inverse transition
matrices respecti\-vely (this implies that transition from  $y^1,\,\ldots,
\,y^n$ to $\ty^1,\,\ldots,\,\ty^n$ is considered as direct transition).
\proclaim{Theorem 1.1}{\sl Under the point transformation \thetag{1.2} the
equations \thetag{1.3} transfer to the equations of the same form
$$
\frac{\partial\ty^i}{\partial\,t}=\sum^n_{j=1}\tilde A^i_j\left(
\frac{\partial^2\ty^j}{\partial x^2}+\shave{\sum^n_{r=1}\sum^n_{s=1}}
\tilde\Gamma^j_{rs}\,\frac{\partial\ty^r}{\partial x}\,
\frac{\partial\ty^s}{\partial x}\right).\hskip -2em
\tag1.4
$$
Parameters $A^i_j$, $\Gamma^j_{rs}$, $\tilde A^i_j$ and $\tilde
\Gamma^j_{rs}$ in \thetag{1.3} and \thetag{1.4} are related as follows
$$
\align
&A^k_i=\sum^n_{m=1}\sum^n_{p=1}S^k_m\,T^p_i\,\tilde A^m_p,
\hskip -2em
\tag1.5\\
&\Gamma^k_{ij}=\sum^n_{m=1}\sum^n_{p=1}\sum^n_{q=1}S^k_m\,T^p_i\,
T^q_j\,\tilde\Gamma^m_{pq}+\sum^n_{m=1}S^k_m\,\frac{\partial T^m_i}
{\partial y^j}.\hskip -2em
\tag1.6
\endalign
$$}
\endproclaim
     If we suppose $y^1,\ldots,y^n$ to be local coordinates on
some manifold $M$, then due to theorem~1.1 parameters $A^k_i$
define a tensor field $\bold A$ of the type $(1,1)$ on $M$.
Parameters $\Gamma^k_{ij}$ in turn can be interpreted as 
components of some affine connection $\Gamma$ on $M$. This
enables us to apply various geometrical methods to the
investigation of the equations \thetag{1.3}. Here we consider
only one problem due to equations \thetag{1.3}, which appears
to have graceful geometrical solution. Similar problems for
other classes of differential equations were considered in
\cite{1--6}.
\head
2. Reduction to the equations of diffusion type.
\endhead
      Equations \thetag{1.1} arise in various models describing
the diffusion phenomena in multicomponent mixtures (see \cite{7--24})
and in theory of integrability (see \cite{25--30}). They are called
the equations {\it of diffusion type}. As for the equations \thetag{1.3} 
those of them, that can be brought to the form \thetag{1.1} by means of 
point transformation \thetag{1.2}, should obviously be involved to the 
class of equations of diffusion type. This gives rise to the problem of 
finding effective description for this class of equations.
\proclaim{Problem}{\sl What condition for tensor field $\bold A$ and
affine connection $\Gamma$ should be fulfilled on $M$ in order to
guarantee the existence of point transformation \thetag{1.2} that 
brings equations \thetag{1.2} to the diffusion form \thetag{1.1}?}
\endproclaim
     Suppose that one could bring \thetag{1.3} to the form \thetag{1.1}
in  local  variables  $\ty^1,\,\ldots,\,\ty^n$. Then in such variables 
for components of affine connection $\Gamma$ we have
$$
\tilde\Gamma^m_{pq}=\frac{1}{2}\sum^n_{s=1}\tilde B^m_s\left(
\frac{\partial\tilde A^s_p}{\partial\ty^q}+\frac{\partial\tilde
A^s_q}{\partial\ty^p}\right).
\tag2.1
$$
Here $\tilde B^m_s$ are the components of the matrix $\tilde B$
which is inverse to the matrix $\tilde A$ formed by components of 
operator field $\bold A$ in local coordinates  $\ty^1,\,\ldots,\,
\ty^n$. Now let's apply formulas \thetag{1.5} in order to transform
the components of the field $\bold A$ in \thetag{2.1}:
$$
\tilde A^s_p=\sum^n_{r=1}\sum^n_{i=1}T^s_r\,S^i_p\,A^r_i,
$$
Further differentiate this relationship with respect to $\ty^q$:
$$
\frac{\partial\tilde A^s_p}{\partial\ty^q}=\sum^n_{r=1}\sum^n_{i=1}
\left(\frac{\partial T^s_r}{\partial\ty^q}\,S^i_p\,A^r_i
+T^s_r\,\frac{\partial S^i_p}{\partial\ty^q}\,A^r_i+T^s_r\,S^i_p
\,\frac{\partial A^r_i}{\partial\ty^q}\right).\hskip -2em
\tag2.2
$$
Let's express the derivatives $\partial A^r_i/\partial\ty^q$, 
$\partial S^i_p/\partial\ty^q$, and $\partial T^s_r/\partial\ty^q$
through partial derivatives in initial variables $y^1,\,\ldots,\,y^n$.
We shall do this by means of formulas
$$
\xalignat 3
&\frac{\partial A^r_i}{\partial\ty^q}=\sum^n_{j=1}
S^j_q\,\frac{\partial A^r_i}{\partial y^j},
&&\ \frac{\partial T^s_r}{\partial\ty^q}=\sum^n_{j=1}
S^j_q\,\frac{\partial T^s_r}{\partial y^j},
&&\ \frac{\partial S^i_p}{\partial\ty^q}=-\sum^n_{j=1}
\sum^n_{\alpha=1}\sum^n_{\beta=1}S^i_\beta\,S^j_q\,
\frac{\partial T^\beta_\alpha}{\partial y^j}\,S^\alpha_p.
\endxalignat
$$
Substitute them into \thetag{2.2}, then substitute the resulting
expression into the formula \thetag{2.1} for the components of
$\Gamma$. Missing expression for the derivative $\partial\tilde A^s_q
/\partial\ty^p$ we get by transposing indices $p$ and $q$:
$$
\align
\tilde\Gamma^m_{pq}=\frac{1}{2}&\sum^n_{r=1}\sum^n_{s=1}\sum^n_{i=1}
\sum^n_{j=1}\tilde B^m_s\,\frac{\partial T^s_r}{\partial y^j}\,
A^r_i\,(S^i_p\,S^j_q+S^i_q\,S^j_p)\,+\\
\vspace{1ex}
+\,\frac{1}{2}&\sum^n_{r=1}\sum^n_{s=1}\sum^n_{i=1}
\sum^n_{j=1}\tilde B^m_s\,T^s_r\,\frac{\partial A^r_i}{\partial y^j}
\,(S^i_p\,S^j_q+S^i_q\,S^j_p)\,-\\
\vspace{1ex}
&-\frac{1}{2}\sum^n_{s=1}\sum^n_{j=1}\sum^n_{\alpha=1}\sum^n_{\beta=1}
\tilde B^m_s\,\tilde A^s_\beta\,\frac{\partial T^\beta_\alpha}
{\partial y^j}\,(S^\alpha_p\,S^j_q+S^\alpha_q\,S^j_p).
\endalign
$$
Now substitute the above expression for $\tilde\Gamma^m_{pq}$ into
the formula \thetag{1.6}. As a result we obtain the following formula 
for $\Gamma^k_{ij}$:
$$
\Gamma^k_{ij}=\sum^n_{r=1}\frac{B^k_r}{2}\left(\frac{\partial A^r_i}
{\partial y^j}+\frac{\partial A^r_j}{\partial y^i}\right)+
\sum^n_{r=1}\sum^n_{p=1}\sum^n_{q=1}B^k_r\,\theta^r_{pq}\,
\frac{A^p_i\,\delta^q_j+A^p_j\,\delta^q_i}{2}.\hskip -2em
\tag2.3
$$
Here $\delta^q_i$ is Kronecker's delta-symbol (representing
components of unit matrix):
$$
\delta^q_i=\cases 1&\text{for \ }q=i,\\
0&\text{for \ }q\neq i.\endcases
$$
By $B^k_r$ we denote the components of operator field $\bold B=
\bold A^{-1}$. The entrance of matrices $S$ and $T$ into \thetag{2.3}
is completely defined by the quantities $\theta^r_{pq}$:
$$
\theta^r_{pq}=\sum^n_{m=1}S^r_m\,\frac{\partial T^m_p}
{\partial y^q}=\sum^n_{m=1}S^r_m\,\frac{\partial^2\ty^m}
{\partial y^p\,\partial y^q}.\hskip -2em
\tag2.4
$$
From \thetag{2.4} we see that $\theta^r_{pq}=\theta^r_{qp}$.
If we know that system of equations \thetag{1.3} admits the
coordinates $\ty^1,\,\ldots,\,\ty^n$, where it takes diffusion 
form \thetag{1.1}, and if these coordinates are given, then
by means of formula \thetag{2.4} we derive $\theta^r_{pq}$ for 
further substitution into the relationship \thetag{2.3}.\par
     Otherwise, when coordinates $\ty^1,\,\ldots,\,\ty^n$ are not
given, one should determine $\theta^r_{pq}$ from \thetag{2.3} and
should consider \thetag{2.4} as the equations that determine
required change of variables \thetag{1.2}. We write these 
equations as 
$$
\frac{\partial T^m_p}{\partial y^q}=\sum^n_{r=1}
\theta^r_{pq}\,T^m_r.\hskip -2em
\tag2.5
$$
Equations \thetag{2.5} form complete system of linear Pfaff
equations respective to the components of the matrix $T$.
They are well known in differential geometry. Compati\-bility
condition for \thetag{2.5} is derived from permutability of 
partial derivatives:
$$
\frac{\partial^2 T^m_p}{\partial y^q\,\partial y^k}=
\frac{\partial^2 T^m_p}{\partial y^k\,\partial y^q}.
\hskip -2em
\tag2.6
$$
Let's calculate both sides of \thetag{2.6} according to equations 
\thetag{2.5}. As a result we get the relationship, which is known
as the condition of ``zero curvature''\footnote{Note that  
``zero curvature" condition \thetag{2.7} does not mean that the
curvature tensor of the affine connection $\Gamma$ is zero.}:
$$
\frac{\partial\theta^m_{kp}}{\partial y^q}-\frac{\partial\theta^m_{qp}}
{\partial y^k}+\sum^n_{r=1}\theta^m_{qr}\,\theta^r_{kp}-\sum^n_{r=1}
\theta^m_{kr}\,\theta^r_{qp}=0.\hskip -2em
\tag2.7
$$
On the base of theory of Pfaff equations (see, for instance, appendix A
in \cite{4}) one can prove the following theorem.
\proclaim{Theorem 2.1}{\sl Let $v$ be some point of the manifold $M$. Matrix
Pfaff equations \thetag{2.5} have solution with non-zero determinant
$\det T\neq 0$ in some neighborhood of the point $v$ on $M$ if and only
if the condition \thetag{2.7} is fulfilled.}
\endproclaim
     Quantities $\theta^r_{pq}$ are symmetric in lower pair of indices
$p$ and $q$. Hence for arbitrary solution of the equations \thetag{2.5} 
we have 
$$
\frac{\partial T^m_p}{\partial y^q}=\frac{\partial T^m_q}{\partial y^p}.
\hskip -2em
\tag2.8
$$
This condition \thetag{2.8} in turn is the compatibility condition
for the equations
$$
\frac{\partial\ty^m}{\partial y^p}=T^m_p(y^1,\ldots,y^n)
$$
which determine new variables $\ty^1,\,\ldots,\,\ty^n$. We shall
formulate the result of the above speculations in form of two
theorems.
\proclaim{Theorem 2.2}{\sl Quantities $\theta^r_{pq}=\theta^r_{qp}$
are defined by some change of local coordinates in some neighborhood 
of the point $v$ on $M$ according to the formulas \thetag{2.4} if
and only if the conditions of ``zero curvature'' \thetag{2.7} are
fulfilled.}
\endproclaim
\proclaim{Theorem 2.3}{\sl System of equations \thetag{1.3} can
be brought to the diffusion form \thetag{1.1} by some change of
local coordinates in some neighborhood of the point $v$ on $M$
if and only if one can find quantities $\theta^r_{pq}=\theta^r_{qp}$
satisfying both conditions \thetag{2.3} and \thetag{2.7} 
simultaneously.}
\endproclaim
     Theorem~2.3 distinguishes two cases in the theory of quasilinear
equations of the form \thetag{1.3}: first is {\bf the case of general 
position} and second is {\bf the degenerate case}. In the case of general 
position quantities $\theta^r_{pq}$ are uniquely determined by the
relationships \thetag{2.3}. Substituting them into \thetag{2.7}
further we can get required equations for $\bold A$ and $\Gamma$
which guarantee that the equations \thetag{1.3} are reducible to 
the form \thetag{1.1}.\par
     In degenerate case relationships \thetag{2.3} do not determine
the quantities $\theta^r_{pq}$ at all or determine with some extent 
of uncertainty (not uniquely). In this case we should analyze 
the relationships \thetag{2.3} for their solvability with respect to 
$\theta^r_{pq}$ and for their compatibility with \thetag{2.7}.
\head
3. Non-degeneracy condition.
\endhead
     The non-degeneracy condition distinguishing the case of general
position, as formulated above, is the condition of unique solvability
of the equations \thetag{2.3} with respect to $\theta^r_{pq}$. These
equations can be written as follows:
$$
\sum^n_{p=1}\sum^n_{q=1}\theta^r_{pq}\,\frac{A^p_i\,\delta^q_j
+A^p_j\,\delta^q_i}{2}=\sum^n_{k=1}A^r_k\,\Gamma^k_{ij}-
\frac{1}{2}\left(\frac{\partial A^r_i}{\partial y^j}
+\frac{\partial A^r_j}{\partial y^i}\right).\hskip -2em
\tag3.1
$$
The equations \thetag{3.1} are linear equations respective to
$\theta^r_{pq}$. Note that they breaks into $n$ separate parts.
For each $r=1,\,\ldots,\,n$ we have closed system of $n^2$ linear 
equations. Matrices of all these systems are the same, they are 
defined by the components of tensor field $\boldsymbol\Lambda$
in the following form:
$$
\Lambda^{pq}_{ij}=\frac{A^p_i\,\delta^q_j+A^p_j\,\delta^q_i
+A^q_i\,\delta^p_j+A^q_j\,\delta^p_i}{4}.
\hskip -2em
\tag3.2
$$
Denote by $w^r_{ij}$ right hand sides of \thetag{3.1}. Then
these equations can be written as 
$$
\sum^n_{p=1}\sum^n_{q=1}\Lambda^{pq}_{ij}\,\theta^r_{pq}=w^r_{ij}.
\hskip -2em
\tag3.3
$$\par
     Let $v$ be some fixed point on the manifold $M$. Denote by $V$
the tangent space to $M$ at this point: $V=T_v(M)$. By $V^*$ denote
the dual space for $V$. Then let's consider a tensor product of two
samples of dual space:
$$
W=V^*\otimes V^*=T^0_2(v,M).
$$
According to this definition $W$ is a space of bilinear forms on the
space $V$. It is represented as a direct sum of two subspaces:
$$
W=W_{\text{sym}}\oplus W_{\text{skew}}.
\hskip -2em
\tag3.4
$$
First is the subspace of symmetric bilinear forms and second is the
subspace of skew-symmetric forms. Tensor $\boldsymbol\Lambda$ can be
treated as a linear operator in the space $W$. Both subspaces in
the sum \thetag{3.4} are invariant with respect to the action of
the operator $\boldsymbol\Lambda\!:W\to W$, restriction of 
$\boldsymbol\Lambda$ to $W_{\text{skew}}$ being identically zero
(since $\Lambda^{pq}_{ij}$ is symmetric in both pairs of indices
upper and lower).\par
     Right hand side of \thetag{3.3} is symmetric in $i$ and $j$,
quantities $\theta^r_{pq}$, that should be determined from the
equations \thetag{3.1}, are also symmetric in their lower indices
$p$ and $q$. Therefore the above condition determining the case of
general position can be formulated as the condition of non-degeneracy
for the restriction of the operator $\boldsymbol\Lambda$ to the
subspace of symmetric bilinear forms:
$$
\boldsymbol\Lambda_{\text{sym}}\!:W_{\text{sym}}\to W_{\text{sym}}.
$$
In other words this is written as
$$
\det\boldsymbol\Lambda_{\text{sym}}\neq 0.\hskip -2em
\tag3.5
$$
Components of the operator $\boldsymbol\Lambda$ are uniquely determined
by the components of $\bold A$ according to the formula \thetag{3.2}.
But nevertheless \thetag{3.5} doesn't follow immediately from $\det\bold 
A\neq 0$. One should check it up separately.\par
\head
4. Effectivization of the non-degeneracy condition.
\endhead
     Let's study the question of how to check up the condition 
\thetag{3.5}, if the matrix of operator field $\bold A$ is given. 
We shall start with the case when operator $\bold A$ has purely 
simple spectrum, i\.~e\. when its eigenvalues $\lambda_1,\,\ldots,
\,\lambda_n$ are distinct and its characteristic polynomial has 
no multiple roots: 
$$
f(\lambda)=\det(\bold A-\lambda\,\bold I)=(-\lambda)^n+\sum^n_{s=1}
\sigma_s(n)\,(-\lambda)^{n-s}=\prod^n_{s=1}(\lambda_s-\lambda).
\hskip -2em
\tag4.1
$$
Here $\sigma_1(n),\,\ldots,\,\sigma_n(n)$ are basic symmetric 
polynomials defined by formula
$$
\sigma_s(n)=\msum_{1\leqslant i_1<\,\ldots\,<i_s\leqslant n}
\lambda_{i_1}\cdot\ldots\cdot\lambda_{i_s}.\hskip -2em
\tag4.2
$$
Last of these polynomials is equal to the determinant of the operator
$\bold A$ (see more details in \cite{31}).\par
     Operator $\bold A$ with simple spectrum in finite dimensional
space $V=T_v(M)$ is diagonalizable, it has a base of eigenvectors
$\bold e_1,\,\ldots,\,\bold e_n$. However, the construction of such 
base may require complexification of the space $V$, since roots
or characteristic polynomial \thetag{4.1} are complex numbers in
general. Denote by $\bold h^1,\,\ldots,\,\bold h^n$ the dual base
in dual space $V^*$. We have the following relationships:
$$
\xalignat 3
&\bold A\bold e_i=\lambda_i\,\bold e_i, &&\bold A^*\bold h^i
=\lambda_i\,\bold h^i,&&\bold h^i(\bold e_k)=\delta^i_k.
\hskip -2em
\tag4.3
\endxalignat
$$
Here $\bold A^*\!:V^*\to V^*$ is a conjugate operator for $\bold A$.
By means of covectors $\bold h^1,\,\ldots,\,\bold h^n$ we shall
build the base in the space of symmetric bilinear forms 
$W_{\text{sym}}$:
$$
\bold w^{ij}=\bold h^i\otimes\bold h^j+\bold h^j\otimes\bold h^i
\text{, \ where \ }1\leqslant i\leqslant j\leqslant n.\hskip -2em
\tag4.4
$$
Denote by $N$ the dimension of the space $W_{\text{sym}}$. On counting 
the number of vectors in the base \thetag{4.4} we find
$$
N=\dim W_{\text{sym}}=\frac{n(n+1)}{2}.\hskip -2em
\tag4.5
$$
Due to \thetag{3.2} and \thetag{4.3} the result of applying
$\boldsymbol\Lambda$ to the vector $\bold h^p\otimes\bold h^q$
is given by
$$
\boldsymbol\Lambda(\bold h^p\otimes\bold h^q)
=\frac{\lambda_p+\lambda_q}{2}\,\bold h^p\otimes\bold h^q.
\hskip -2em
\tag4.6
$$
Now let's apply the operator \thetag{4.6} to the vectors of
the base \thetag{4.4}. By means of direct calculations we get
$$
\boldsymbol\Lambda_{\text{sym}}(\bold w^{ij})=\boldsymbol\Lambda
(\bold w^{ij})=\frac{\lambda_i+\lambda_j}{2}\,\bold w^{ij}.
\hskip -2em
\tag4.7
$$
The relationship \thetag{4.7} means that operator $\boldsymbol
\Lambda_{\text{sym}}$ is diagonalizable, base \thetag{4.4}
being the base of eigenvectors of this operator. Eigenvalues
$\mu_{ij}$ corresponding to the eigenvectors $\bold w^{ij}$
are given by the formula
$$
\mu_{ij}=\frac{\lambda_i+\lambda_j}{2}\text{, \ where \ }
1\leqslant i\leqslant j\leqslant n.\hskip -2em
\tag4.8
$$
Using \thetag{4.8} we can calculate the determinant of the operator
$\boldsymbol\Lambda_{\text{sym}}$:
$$
\det\boldsymbol\Lambda_{\text{sym}}=\sigma_N(N)=\prod_{i\leqslant
j}\left(\frac{\lambda_i+\lambda_j}{2}\right).\hskip -2em
\tag4.9
$$\par
     Let's introduce lexicographic ordering in the set of pairs
of indices $(i,j)$, where $i\leqslant j$. Say that $(i_1,j_1)<(i_2,j_2)$
if $j_1<j_2$ or if $i_1<i_2$ in the case of coincidence $j_1=j_2$. Then
consider the sum
$$
\sigma_s(N)=\msum_{(i_1,j_1)<\,\ldots\,<(i_s,j_s)}\frac{(\lambda_{i_1}
+\lambda_{j_1})\cdot\ldots\cdot(\lambda_{i_s}+\lambda_{j_s})}{2^s}.
\hskip -2em
\tag4.10
$$
Formula \thetag{4.10} determines basic symmetric polynomials for
the set of numbers $\mu_{ij}$ from \thetag{4.8}. For $s=N$ it is
equivalent to \thetag{4.9}. Its easy to note, that right hand sides
of \thetag{4.10} do not change under any transposition of numbers
$\lambda_1,\,\ldots,\,\lambda_n$. They are symmetric polynomials
in $\lambda_1,\,\ldots,\,\lambda_n$. Hence according the
well known theorem (see \cite{31}, chapter~6, \S\,2) one can express 
them trough basic symmetric polynomials:
$$
\sigma_s(N)=P_s(\sigma_1(n),\ldots,\sigma_n(n))\text{, \ where \ }
s=1,\,\ldots,\,N.\hskip -2em
\tag4.11
$$
For any particular $n$ we can derive explicit expressions for 
polynomials $P_s(\sigma_1,\ldots,\sigma_n)$. For instance, if
$n=2$, we have
$$
\pagebreak 
\xalignat 3
&\quad P_1=\frac{3}{2}\,\sigma_1,
&&P_2=\frac{1}{2}\,\sigma_1^2+\sigma_2,
&&P_3=\frac{1}{2}\,\sigma_2\,\sigma_1,
\hskip -3em
\tag4.12
\endxalignat
$$
From \thetag{4.5} here we get $N=3$. Since $\sigma_2(2)=
\det A$ is non-zero, for $n=2$ the non-degeneracy condition 
\thetag{3.5} can be written as
$$
\sigma_1(2)=\tr\bold A\neq 0.\hskip -2em
\tag4.13
$$
\proclaim{Theorem 4.1}{\sl For $n=2$ the system of equations
\thetag{1.3} belongs to the non-degenerate class if and only if
$\det\bold A\neq 0$ and $\tr\bold A\neq 0$.}
\endproclaim
     The relationship \thetag{4.13} means that operator $\bold A$
shouldn't have two eigenvalues which differ only in sign. In this
form the non-degeneracy condition can be stated in multidimensional
case as well. For arbitrary $n$ it follows from \thetag{4.9}.
Indeed, we have 
$$
\prod_{i\leqslant j}\left(\frac{\lambda_i+\lambda_j}{2}\right)=
\det\bold A\cdot\prod_{i<j}\left(\frac{\lambda_i+\lambda_j}{2}
\right)\neq 0,\hskip -3em
\tag4.14
$$
which is equivalent to $\det\bold A\neq 0$ and $\lambda_i\neq-\lambda_j$
for $i\neq j$.\par
     In order to check up the condition \thetag{4.14} one need not
find eigenvalues of the operator $\bold A$. This condition can be
written in terms of coefficients of characteristic polynomial
\thetag{4.1}, provided we have explicit expressions for polynomials
in \thetag{4.11}. The latter is algorithmically solvable problem.\par
     Another effective way of checking up the condition \thetag{4.14}
is based on the theory of resultants. Let $f(\lambda)$ be a 
characteristic polynomial for the operator $\bold A$. Its roots are
non-zero since $\det\bold A\neq 0$, therefore the condition
$\lambda_i\neq-\lambda_j$ for $i\neq j$ means that polynomials
$f(\lambda)$  and  $f(-\lambda)$  have  no  common  roots.   Then 
resultant
of these two polynomials is non-zero (see \cite{31} or \cite{32}):
$$
\Res_\lambda\bigl[\det(\bold A-\lambda\,\bold I),\det(\bold A
+\lambda\,\bold I)\bigr]\neq 0.\hskip -3em
\tag4.15
$$
\proclaim{Theorem 4.2}{\sl System of equations \thetag{1.3} belongs
to the non-degenerate class if and only if $\det\bold A\neq 0$ and
condition \thetag{4.15} is fulfilled.}
\endproclaim
\head
5. Finding inverse operator for $\boldsymbol\Lambda_{\text{sym}}$.
\endhead
     Suppose that non-degeneracy condition stated in theorem~4.2
is fulfilled. Let's find regular way for inverting the operator
$\boldsymbol\Lambda_{\text{sym}}\!:W_{\text{sym}}\to W_{\text{sym}}$. 
One need to do it for to find $\theta^r_{pq}$ from the equations
\thetag{3.3}. Let $f(\lambda)$ be a characteristic polynomial for
the operator $\bold A$. Lets consider pair of polynomials
$f(\lambda+\mu)$ and $f(\mu-\lambda)$, where $\mu$ is a parameter.
These polynomials are factored as follows:
$$
\aligned
f(\lambda+\mu)&=(-1)^n\prod^n_{i=1}\bigl(\lambda-(\lambda_i-\mu)
\bigr)\\
f(\mu-\lambda)&=\prod^n_{i=1}\bigl(\lambda-(\mu-\lambda_i)\bigr).
\endaligned\hskip -2em
\tag5.1
$$
For to calculate the resultant of polynomials \thetag{5.1} let's
apply the following theorem of algebra (see \cite{31} or \cite{32}
for proof).
\proclaim{Theorem 5.1}{\sl Suppose that two polynomials $g(\lambda)$ 
and $h(\lambda)$  with  leading  coefficients  $a$  and  $b$  are 
factored
to a product of linear terms
$$
\aligned
g(\lambda)&=a\cdot(\lambda-\alpha_1)\cdot\ldots\cdot(\lambda
-\alpha_n),\\
h(\lambda)&=b\cdot(\lambda-\beta_1)\cdot\ldots\cdot(\lambda
-\beta_m).
\endaligned\hskip -2em
\tag5.2
$$
Then resultant of these polynomials is given by the formula
$$
\Res_\lambda\bigl[g(\lambda),h(\lambda)\bigr]=a^n\,b^m\,
\prod^n_{i=1}\prod^m_{j=1}(\alpha_i-\beta_j).
\hskip -2em
\tag5.3
$$}
\endproclaim
\noindent
Comparing \thetag{5.1} with \thetag{5.2} from \thetag{5.3}
we obtain the resultant of polynomials \thetag{5.1}:
$$
\gathered
\Res_\lambda\bigl[f(\mu-\lambda),f(\lambda+\mu)\bigr]=
(-1)^n\prod^n_{i=1}\prod^n_{j=1}(2\,\mu-\lambda_i-\lambda_j)=\\
\vspace{1ex}
=2^{n^2}\,\prod_{i<j}\left(\frac{\lambda_i+\lambda_j}{2}
-\mu\right)\,\prod_{i\leqslant j}\left(\frac{\lambda_i+\lambda_j}
{2}-\mu\right).
\endgathered\hskip -2em
\tag5.4
$$
Let's compare \thetag{5.4} with formula \thetag{4.8} for eigenvalues
of the operator $\boldsymbol\Lambda_{\text{sym}}$. This comparison
gives the formula 
$$
\varphi(\mu)^2=2^{-n^2}\,f(\mu)\,\Res_\lambda\bigl[f(\mu-\lambda),
f(\lambda+\mu)\bigr].\hskip -2em
\tag5.5
$$
Here $\varphi(\mu)=\det(\boldsymbol\Lambda_{\text{sym}}-\mu\,\bold I)$
is a characteristic polynomial for the operator $\boldsymbol
\Lambda_{\text{sym}}$. Its remarkable that formula \thetag{5.5} express
the square of characteristic polynomial of operator $\boldsymbol 
\Lambda_{\text{sym}}$ in terms of characteristic polynomial of operator 
$\bold A$.\par
     Denote by $\varepsilon(\mu)$ the polynomial in right hand side of
\thetag{5.5}. According to the well known theorem of Hamilton and 
Cayley
(see \cite{32} or \cite{33}) when substituting the operator $\boldsymbol
\Lambda_{\text{sym}}$ into its characteristic polynomial $\varphi(\mu)$
we should get zero operator. This is true for the polynomial 
$\varepsilon(\mu)=\varphi(\mu)^2$ as well:
$$
\varepsilon(\boldsymbol\Lambda_{\text{sym}})=\sum^M_{i=0}\varepsilon_i
\,\boldsymbol\Lambda_{\text{sym}}^i=0.\hskip -2em
\tag5.6
$$
Here $M=2\,N=n(n+1)$ is the degree of polynomial $\varepsilon(\mu)=
\varphi(\mu)^2$. Now we can rewrite \thetag{5.6} as follows
$$
\left(\,\shave{\sum^M_{i=1}}\varepsilon_i\,
\boldsymbol\Lambda_{\text{sym}}^{i-1}
\right)\cdot\boldsymbol\Lambda_{\text{sym}}
=-\varepsilon_0\,\bold I.\hskip -2em
\tag5.7
$$
Quantity $\varepsilon_0$ in right hand side of \thetag{5.7} is
the value of polynomial $\varepsilon(\lambda)$ at $\lambda=0$:
$$
\varepsilon_0=\varepsilon(0)=\varphi(0)^2=(\det\boldsymbol
\Lambda_{\text{sym}})^2.
$$
Due to non-degeneracy condition \thetag{3.5} the quantity $\varepsilon_0$ 
is non-zero. Therefore we can use \thetag{5.7} in order to calculate the 
operator $\boldsymbol\Lambda_{\text{sym}}^{-1}$:
$$
\boldsymbol\Lambda_{\text{sym}}^{-1}=-\sum^M_{i=1}\frac{\varepsilon_i}
{\varepsilon_0}\,\boldsymbol\Lambda_{\text{sym}}^{i-1}.
$$
We can expand operator $\boldsymbol\Lambda_{\text{sym}}^{-1}\!:
W_{\text{sym}}\to W_{\text{sym}}$ to be identically zero in the subspace 
of skew-symmetric bilinear forms $W_{\text{skew}}$. Though being 
degenerate, expanded operator $\boldsymbol\Lambda_{\text{sym}}^{-1}$ is 
defined in the whole space of bilinear forms $W$. One can define 
this operator by the following formula
$$
\boldsymbol\Lambda_{\text{sym}}^{-1}=-\sum^M_{i=1}\frac{\varepsilon_i}
{\varepsilon_0}\,\boldsymbol\Lambda^{i-1}.\hskip -2em
\tag5.8
$$
Right hand side of \thetag{5.8} is interpreted as a tensor field
$\bold D$ of the type $(2,2)$ whose components $D^{ij}_{pq}$ are
symmetric in upper and lower indices.
\proclaim{Theorem 5.2}{\sl If $\det A\neq 0$ and if non-degeneracy 
condition \thetag{3.5} written in the form \thetag{4.15} is
fulfilled, then equations \thetag{3.3} can be solved with respect
to $\theta^r_{pq}$:
$$
\theta^r_{pq}=\sum^n_{i=1}\sum^n_{j=1}D^{ij}_{pq}\,w^r_{ij}.\hskip -2em
\tag5.9
$$
Here $D^{ij}_{pq}$ are components of tensor field $\bold D=\boldsymbol
\Lambda_{\text{sym}}^{-1}$ defined by the formula \thetag{5.8}.}
\endproclaim
    Naturally we can choose more straightforward way of calculating
components of tensor field $\bold D$ other than formula \thetag{5.8}.
In order to do it we should only remember that $\bold D$ defines
an operator inverse to the operator $\boldsymbol\Lambda_{\text{sym}}$ 
in subspace $W_{\text{sym}}$:
$$
\bold D\cdot\boldsymbol\Lambda_{\text{sym}}=\bold I.\hskip -2em
\tag5.10
$$
Written in coordinates the relationship \thetag{5.10} looks like
$$
\sum^n_{i=1}\sum^n_{j=1}D^{ij}_{rs}\,\Lambda^{pq}_{ij}=\frac{\delta^p_r\,
\delta^q_s+\delta^p_s\,\delta^q_r}{2}.
$$
This is the system of $N^2$ linear equations with respect to $D^{ij}_{rs}$, 
where $N$ is defined by \thetag{4.5}. Non-degeneracy condition \thetag{4.15} 
guarantees the compatibility of these equations and the uniqueness of their 
common solution.
\head
6. Criterion for being in diffusion class.
\endhead
    On solving the equations \thetag{3.3} respective to $\theta^r_{pq}$
in the form \thetag{5.9} we can substitute this expression into \thetag{2.7} 
and we can apply theorem~2.3. 
\proclaim{Theorem 6.1}{\sl In non-degenerate case equations \thetag{1.3} 
can be brought to the diffusion form \thetag{1.1} by means of some point 
transformation \thetag{1.2} if and only if the quantities
$$
\theta^r_{pq}=\sum^n_{i=1}\sum^n_{j=1}D^{ij}_{rs}\left(\,
\shave{\sum^n_{k=1}}A^r_k\,\Gamma^k_{ij}-
\frac{1}{2}\left(\frac{\partial A^r_i}{\partial y^j}
+\frac{\partial A^r_j}{\partial y^i}\right)\right).\hskip -2em
\tag6.1
$$
do satisfy ``zero curvature'' condition \thetag{2.7}.}
\endproclaim
     Due to formulas derived in sections~3, 4 and 5 we can check up 
non-degeneracy condition and we can explicitly calculate
$D^{ij}_{rs}$ in those variables $y^1,\,\ldots,\,y^n$, where the 
components of operator field $\bold A$ are initially given. Therefore
theorem~2.3 is an effective criterion for testing if the equation
belong to the diffusion class in non-degenerate case. The algorithm
given by this theorem can be easily realized by means of any of now
existing program packages for symbolic calculations.\par
    {\bf Note}. {\sl All results in sections~3, 4 and 5 were obtained
under the assumption that operator $\bold A$ has purely simple
spectrum. However, they remain true for arbitrary operators, since
any operator with degenerate spectrum can be obtained as a limit of
some sequence of operators with purely simple spectrum.}
\head
7. Explicit formulas for $n=2$.
\endhead
     As it was shown above, in two-dimensional case the non-degeneracy
condition reduces to the following relationships
$$
\xalignat 2
&\det\bold A\neq 0,&&\tr\bold A\neq 0
\endxalignat
$$
(see theorem~4.1). Using \thetag{4.12} we can find explicit expression
for the characteristic polynomial $\varphi(\mu)$ of the operator
$\boldsymbol\Lambda_{\text{sym}}$. Its coefficients are expressed
in terms of parameters $\sigma_1=\tr\bold A$ and $\sigma_2=\det\bold A$:
$$
\varphi(\mu)=-\mu^3+\frac{3\,\sigma_1}{2}\,\mu^2-\frac{\sigma_1^2
+2\,\sigma_2}{2}\,\mu+\frac{\sigma_1\,\sigma_2}{2}.\hskip -2em
\tag7.1
$$
From operator equality $\varphi(\boldsymbol\Lambda_{\text{sym}})=0$
due to \thetag{7.1} we get the following expression:
$$
\boldsymbol\Lambda_{\text{sym}}^{-1}=\frac{2\,\boldsymbol\Lambda
^2-3\,\sigma_1\,\boldsymbol\Lambda+(\sigma_1^2+2\,\sigma_2)\,\bold I}
{\sigma_1\,\sigma_2}.\hskip -2em
\tag7.2
$$
Substituting \thetag{3.2} into \thetag{7.2} we obtain the components of
the tensor field $\bold D=\boldsymbol\Lambda_{\text{sym}}^{-1}$:
$$
\xalignat 2
&D^{11}_{11}=\frac{\sigma_2+(A^2_2)^2}{\sigma_1\,\sigma_2},
&&D^{22}_{22}=\frac{\sigma_2+(A^1_1)^2}{\sigma_1\,\sigma_2},\\
\vspace{1ex}
&D^{11}_{22}=\frac{(A^1_2)^2}{\sigma_1\,\sigma_2},
&&D^{22}_{11}=\frac{(A^2_1)^2}{\sigma_1\,\sigma_2},\\
\vspace{1ex}
&D^{12}_{11}=D^{21}_{11}=-\frac{A^2_1\,A^2_2}{\sigma_1\,\sigma_2},
&&D^{12}_{22}=D^{21}_{22}=-\frac{A^1_2\,A^1_1}{\sigma_1\,\sigma_2},\\
\vspace{1ex}
&D^{11}_{12}=D^{11}_{21}=-\frac{A^1_2\,A^2_2}{\sigma_1\,\sigma_2},
&&D^{22}_{12}=D^{22}_{21}=-\frac{A^2_1\,A^1_1}{\sigma_1\,\sigma_2},\\
\allowdisplaybreak
\vspace{1ex}
&D^{12}_{12}=D^{21}_{12}=\frac{A^1_1\,A^2_2}{\sigma_1\,\sigma_2},
&&D^{12}_{21}=D^{21}_{21}=\frac{A^1_1\,A^2_2}{\sigma_1\,\sigma_2}.
\endxalignat
$$
Now we should substitute the above expressions into \thetag{6.1}
for to find parameters $\theta^r_{pq}$ and substitute them into
\thetag{2.7}. As a result we obtain four differential equations
binding $A^i_j$ and $\Gamma^k_{ij}$. If they hold, then the
equations \thetag{1.3} can be brought to the form \thetag{1.1} 
in two-dimensional case $n=2$.


\Refs
\ref\no 1\by Pavlov~M\.~V\., Svinolupov~S\.~I\., Sharipov~R\.~A\.\paper
Invariant criterion of hydrodynamical integrability for the equations
of hydrodynamical type\jour Funk\. analiz i pril\.\yr 1996\vol 30
\issue 1\pages 18--29\moreref\inbook see also in book ``Integrability
in dynamical systems''\publ Inst. of Math\. and IC RAS\publaddr
Ufa\yr 1994\pages 27--48\moreref and see also in Electronic archive 
at LANL, solv-int \#9407003
\endref
\ref\no 2\by Ferapontov~E\.~V\., Sharipov~R\.~A\.\paper On first-order 
conservation laws for systems of hydrodynamic-type equations
\jour TMF\yr 1996\vol 108\issue 1\pages 109--128
\endref
\ref\no 3\by Dmitrieva~V\.~V\., Sharipov~R\.~A\.\paper On the point
transformations for the second order differential equations\inbook
Electronic archive at LANL (1997), solv-int \#9703003\pages 1--14
\endref
\ref\no 4\by Sharipov~R\.~A\.\paper On the point transformations for
the equation $y''=P+3\,Q\,y'+3\,R\,{y'}^2+S\,{y'}^3$\inbook
Electronic archive at LANL (1997), solv-int \#9706003\pages 1--35
\endref
\ref\no 5\by Mikhailov~O\.~N\., Sharipov~R\.~A\.\paper On the point
expansion for certain class of differential equations of second 
order\inbook Electronic archive at LANL (1997), solv-int
\#9712001\pages 1--8
\endref
\ref\no 6\by Sharipov~R.~A.\paper Effective procedure of point
classification for the equations $y''=P+3\,Q\,y'+3\,R\,{y'}^2
+S\,{y'}^3$\inbook Electronic archive at LANL (1998), Math\.DG
\#9802029\pages 1--35
\endref
\ref\no 7\by Hoff~D\., Zumbrun~K\.\paper Multi-dimensional diffusion
waves for Navier-Stokes equations of compressible flow\jour Indiana
Univ\. Math\. Journ\.\vol 44\issue 2\yr 1995\pages 603--676
\endref
\ref\no 8\by Broadbridge~P\., White~I.\yr 1987\paper Constant rate
rainfall infiltration: a versatile nonlinear model. I. Analytic
solution\jour Water Resources Res.\vol 24\pages 145--154
\endref
\ref\no 9\by De~Vries~D\.~A\.\yr 1958\paper Simultaneous transfer of
heat and moisture in porous media\jour Trans\. of American Geophys\.
Union\vol 39\pages 909--916
\endref
\ref\no 10\by Edwards~M\.~P\., Broadbridge~P\.\yr 1994\paper Exact
transient solutions to nonlinear diffusion-convection equations in
higher dimensions \jour Journ. of Phys\.~A: Math\. Gen\.\vol 27
\pages 5455--5465
\endref
\ref\no 11\by Jackson~R\.~D\.\yr 1973\paper Diurnal soil watertime-depth
patterns\jour Proc\. of American Soil Science Soc\.\vol 37\pages 505-509
\endref
\ref\no 12\by Jury~W\.~A\., Letey~J\., Stolzy~L\.H\.\yr 1981\paper Flow
of water and energy under desert conditions\inbook Water in Desert
Ecosystems\eds D\.~D\.~Evans and J\.~L\.~Thames\publ Dowden,
Hutchinson and Ross\publaddr Stroudsburg, PA\pages 92--113
\endref
\ref\no 13\by Knight~J\.~H\., Philip~J\.~R\.\yr 1974\paper Exact solutions
in nonlinear diffusion\jour Journ\. Eng\. Math\.\vol 8\pages 219--227
\endref
\ref\no 14\by Myerscough~M\.~R\., Murray~J\.~D\.\yr 1992\paper Analysis of
propagating pattern in a chemotaxis system\jour Bull\. Math\. Biol\.
\vol 54\issue 1\pages 77--94
\endref
\ref\no 15\by Ovsyannikov~L\.~V\.\yr 1959\paper Group properties of
nonlinear heat equation\jour Dokladi AN SSSR\vol 123\issue 3\pages 
492--495
\endref
\ref\no 16\by Philip~J\.~R\.\yr 1988\paper Quasianalytic and analytic
approaches to unsaturated flow\inbook Flow and Transport in the Natural
Environment: Advances and Applications\eds W\.~L\.~Steffen and
O\.~T\.~Denmead\publ Springer\publaddr Berlin\pages 30--48
\endref
\ref\no 17\by Philip~J\.~R\., De~Vries~D\.~A\.\yr 1957\paper Moisture
movement porous materials under a temperature gradient\jour Trans. of
American Geophys. Union\vol 38\pages 222--232
\endref
\ref\no 18\by Rose~C\.~W\.\yr 1968\paper Water transport in soil with
a daily temperature wave. I.\jour Aust\. Jour\. of Soil Res\.\vol 6
\pages 31--44
\endref
\ref\no 19\by Shepherd~R\., Wiltshire~R\.~J\.\yr 1995\paper An
analytical approach to coupled heat and moisture transport in soil
\jour Transport Porous Media\vol 20\pages 281--304
\endref
\ref\no 20\by Shepherd~R\., Wiltshire~R\.~J\.\yr 1996\paper Spectral
decompositions in nonlinear coupled diffusion\jour IMA Journ\. of Appl\.
Math\.\vol 56\pages 277--287
\endref
\ref\no 21\by Sherratt~J\.~A\.\yr 1994\paper Chemotaxis and chemokinesis
in eukaryotic cells: the Keller-Segel equations as an approximation to
a detailed model\jour Bull\. Math\. Biol\.\vol 56\issue 1\pages 129--146
\endref
\ref\no 22\by Sophocleous~C\.\yr 1996\paper Potential symmetries of
nonlinear diffusion-convection equations\jour Journ\. of Phys\. A:
Math\. Gen\.\vol 29\pages 6951-6959
\endref
\ref\no 23\by Wiltshire~R\.~J\.\yr 1994\paper The use of Lie transformation
groups in the solution of the coupled diffusion equation\jour Journ\. of
Phys\. A: Math\. Gen\.\vol 27\pages 7821-7829
\endref
\ref\no 24\by Baikov~V\.~A\., Gladkov~A\.~V\., Wiltshire~R\.~J\.\yr 1998
\paper Lie symmetry classification analysis for nonlinear coupled 
diffusion\jour Journ\. of Phys\. A: Math\. Gen\.\vol 31\pages 7483-7499
\endref
\ref\no 25\by Mikhailov~A\.~V\., Shabat~A\.~B\.\paper Integrability
conditions for the system of two equations of the form $u_t=A(u)\,u_{xx}
+F(u,u_x)$\inbook Preprint of Landau Inst for Theor\. Phys\.,
Academy of Science of USSR\yr 1985\publaddr Ghernogolovka
\endref
\ref\no 26\by Svinolupov~S\.~I\.\yr 1987\paper On multi-field
analogs of Burgers equation\inbook Preprint of report to Presidium
of Bashkir branch of Academy of Sciences of USSR
\endref
\ref\no 27\by Svinolupov~S\.~I\.\paper On the analogs of the Burgers
equation\jour Phys\. Lett\. A\vol 135\issue 1\yr 1989\page 32
\endref
\ref\no 28\by Habibullin~I\.~T\., Svinolupov~S\.~I\. Integrable boundary
conditions for multi-field Burgers equation\jour Mat\. Zametki
\toappear
\endref
\ref\no 29\by Golubchik~I\.~Z\., Sokolov~V\.~V\.\paper Integrable
equations generated by constant solutions of Yang-Baxter equation
\jour Funk\. analiz i pril\.\vol 30\issue 4\yr 1996\pages 68-71
\endref
\ref\no 30\by Golubchik~I\.~Z\., Sokolov~V\.~V\.\paper Integrable 
equations on the graded Lie algebras\jour TMF\vol 112\issue 3
\yr 1997\pages 375-383
\endref
\ref\no 31\by Kostrikin~A\.~I\.\book Introduction to algebra\yr 1977
\publ Nauka\publaddr Moscow
\endref
\ref\no 32\by Lang~S\.\book Algebra\yr 1965\publ Addison-Wesley
Publishing Co\.\publaddr New York
\endref
\ref\no 33\by Sharipov~R\.~A\.\book Course of linear algebra and 
multidimensional geometry\yr 1996\publ Bashkir State University
\publaddr Ufa
\endref
\endRefs
\enddocument
\end